\documentclass[12pt,reqno]{article}
\usepackage{amssymb}
\usepackage{amsthm}
\usepackage{amsmath}
\newtheorem{theorem}{Theorem}[section]
\newtheorem{corollary}[theorem]{Corollary}

\newtheorem{examples}[theorem]{Examples}
\newtheorem{remark}[theorem]{Remark}

\begin{document}
\begin{center} \Large{\bf $\lambda$-symmetry criteria for linearization  of second order ODEs via point transformations }
\end{center}
\medskip

\begin{center}
Ahmad Y. Al-Dweik$^*$, M. T. Mustafa$^*$, Raed A. Mara'beh$^*$ and F. M. Mahomed$^{**}$\\

{$^*$Department of Mathematics \& Statistics, King Fahd University
of Petroleum and Minerals, Dhahran 31261, Saudi Arabia}\\
{$^{**}$School of Computational and Applied Mathematics, Centre for
Differential Equations, Continuum Mechanics and Applications,
University of the Witwatersrand, Wits 2050,  South Africa}\\

aydweik@kfupm.edu.sa, tmustafa@kfupm.edu.sa,
raedmaraabeh@kfupm.edu.sa and Fazal.Mahomed@wits.ac.za
\end{center}
\begin{abstract}
An alternative proof of Lie's approach for linearization of scalar second
order ODEs is derived using the relationship between
$\lambda$-symmetries and first integrals. This relation
further leads to a new $\lambda$-symmetry linearization criteria
for second order ODEs which provides a new approach for
constructing the linearization transformations with lower
complexity. The effectiveness of the approach is illustrated by
obtaining the local linearization transformations for the
linearizable nonlinear ODEs of the form $y''+F_1(x,y)y'+F(x,y)=0$.
Examples of linearizing nonlinear ODEs which are quadratic or
cubic in the first derivative are also presented.
\end{abstract}
\section{Introduction}
The initial seminal studies of scalar second-order
ordinary differential equations (ODEs) which are linearizable by means of point transformations are due to Lie \cite{Lie} and Tress\'e \cite{Tre}.
In recent decades there have been a resurgence of interest in this topic (see the reviews \cite{IMah, Mah, Qad}).

It was shown by Lie  \cite{Lie} that any second-order ODE
\begin{equation}\label{a}
y''=f(x,y,y')
\end{equation}
which is linearizable via point transformations is at most cubic
in the first derivative, i.e. it has the form
\begin{equation}\label{a0}
y''+F_3(x,y){y'}^3+F_2(x,y){y'}^2+F_1(x,y)y'+F(x,y)=0.
\end{equation}
It is well known \cite{Nail} that any second order linear ODE can
be transformed via point transformations to the free particle
equation
\begin{equation}\label{a1}
\begin{array}{l}
 u_{tt}=0.\\
\end{array}
\end{equation}
Therefore, all linearaizable second-order ODEs  (\ref{a}) are
obtained from the free particle equation (\ref{a1}) via point
transformations. Precisely, the free particle equation (\ref{a1})
can be transformed by an arbitrary change of variables
\begin{equation}\label{a2}
t=\phi \left( x,y  \right),~u =\psi \left( x,y  \right),~~~J=\phi_x\psi_y   -  \phi_y \psi_x\neq0,\\
\end{equation}
where $J$ is the Jacobian, to the family of ODEs (\ref{a0}) with
the coefficients $F(x,y), F_1(x,y),$ $F_2(x,y)$ and $F_3(x,y)$
satisfying the following system of partial differential equations
\begin{equation}\label{a3}
\begin{array}{llll}
F_3(x,y)=A,&F_2(x,y)=B+2w,& F_1(x,y)=P+2z,&F(x,y)=Q,
\end{array}
\end{equation}
in which
\begin{equation}\label{a4}
\begin{large}
\begin{array}{lll}
\mbox{\normalsize{\textit{A}=}}\frac { \phi_y \psi_{yy}   - \psi_y  \phi_{yy}  }{ \phi_x\psi_y   -  \phi_y  \psi_x  },& \mbox{\normalsize{\textit{B}=}}\frac { \phi_x \psi_{yy}   - \psi_x  \phi_{yy}  }{ \phi_x\psi_y   -  \phi_y  \psi_x  },&\mbox{\normalsize{\textit{w}=}}\frac { \phi_y \psi_{xy}   - \psi_y  \phi_{xy}  }{ \phi_x\psi_y   -  \phi_y  \psi_x  },\\
\mbox{\normalsize{\textit{Q}=}}\frac { \phi_x \psi_{xx}   - \psi_x  \phi_{xx}  }{ \phi_x\psi_y   -  \phi_y  \psi_x  },&\mbox{\normalsize{\textit{z}=}}\frac { \phi_x \psi_{xy}   - \psi_x  \phi_{xy}  }{ \phi_x\psi_y   -  \phi_y  \psi_x  },& \mbox{\normalsize{\textit{P}=}}\frac { \phi_y \psi_{xx}   - \psi_y  \phi_{xx}  }{ \phi_x\psi_y   -  \phi_y  \psi_x  }. \\
\end{array}
\end{large}
\end{equation}
Lie \cite{Lie} crucially also found the following
over-determined system of four equations
\begin{equation}\label{b1}
\begin{array}{lll}
w_x=\,z  w -F \,{\it F_3}    -\frac{1}{3}{\frac {\partial F_1}{\partial y}} +\frac{2}{3}{\frac {\partial F_2}{\partial x _{{}}}},\\
w_y=-  w^{2} +  {\it F_2} w +  {\it F_3}z+ {\frac {\partial F_3}{\partial x}}  -{\it F_1} {\it F_3},  \\
z_x=z^{2} -  {\it F_1} z -  F w  + {\frac {\partial F}{\partial y_{{}}}} +F {\it F_2},  \\
z_y=-\,z  w +F \,{\it F_3}    -\frac{1}{3}{\frac {\partial F_2}{\partial x}} +\frac{2}{3}{\frac {\partial F_1}{\partial y _{{}}}},\\
\end{array}
\end{equation}
which are called the Lie conditions. The compatibility of Lie's
conditions give the following well known Lie-Tress\'e linearization test
for ODEs of the form (\ref{a0}), viz.
\begin{equation}\label{b2}
\begin{array}{lll}
{\frac {\partial ^{2} F_1}{\partial {y}^{2}}} -2\,{\frac {\partial ^{2} F_2}{\partial y\partial x}}  +3\,{\frac {\partial ^{2}F_3}{\partial{x}^{2} }}-3\,  {\frac {\partial F_1}{\partial x}}   {\it F_3}  -3\,{\it F_1}  {\frac {\partial F_3}{\partial x}} +6\,  {\frac {\partial F}{\partial y}}   {\it F_3}  +3\,F  {\frac {\partial F_3}{\partial y}}  -{\it F_2} {\frac {\partial F_1}{\partial y}}  + 2\,{\it F_2}{\frac {\partial F_2}{\partial x}}=0  \\
{\frac {\partial ^{2} F_2}{\partial {x}^{2}}} -2\,{\frac {\partial ^{2} F_1}{\partial y\partial x}}  +3\,{\frac {\partial ^{2}F}{\partial{y}^{2} }}+3\,  {\frac {\partial F}{\partial y}}   {\it F_2}  +3\,{\it F}  {\frac {\partial F_2}{\partial y}} -3\,  {\frac {\partial F}{\partial x}}   {\it F_3}  -6\,F  {\frac {\partial F_3}{\partial x}}  +{\it F_1} {\frac {\partial F_2}{\partial x}}  - 2\,{\it F_1}{\frac {\partial F_1}{\partial y}}=0.  \\
\end{array}
\end{equation}
It was Tress\'e \cite{Tre} who first obtained the invariant criteria (\ref{b2}).

We point out that the conditions (\ref{b2})  can also be deduced
via the Cartan equivalence approach \cite{Gri}. More recently in
two independent papers \cite{Ibr, MahQ}, the authors present
geometrical proofs for the determination of the invariant
conditions (\ref{b2}). It should be noted that the projections
utilized in \cite{Ibr, MahQ} were different to each other. In
\cite{Ali} Lie's linearizability criteria were utilized to
determine  a linearizable class of a system of two second order equations that are
obtainable from complex scalar second order ODEs. Furthermore,
Mahomed and Qadir \cite{MahQ1} found criteria for conditional
linearizability  of third order ODEs via point transformation
subject to a Lie-Tress\'e linearizable second order ODE.

One can mention here as well the results on the algebraic criteria
for linearization by point transformations of scalar second order
ODEs (\ref{a}). It is well known from Lie's work that such
linearizable ODEs possess eight point symmetries. The question
arises if one can conclude linearization of (\ref{a}) in case one
has knowledge of fewer than eight symmetries of the ODE (\ref{a}).
The answer is affirmative. In fact Lie was the first to obtain
algebraic criteria for linearization when the ODE (\ref{a}) admits
two symmetries which are connected. Some further studies
\cite{Sar, Lea} provided input on the algebraic criteria for
linearization when the ODE (\ref{a}) has two unconnected
symmetries.

Lie proved that the ODE (\ref{a}) is linearizable by point
transformations if and only if the over-determined system
(\ref{b1}) is compatible \cite{Lie}. In order to prove the
sufficiency of the compatibility of (\ref{b1}) for linearization,
Lie showed that the system (\ref{b1}) can be linearized and the
resulting linear system belongs to a class of special type of
linear systems which can be reduced to a linear third-order ODE.
Thus the three solutions $(z, w), (z_1, w_1)$ and $(z_2,w_2)$ of
the system (\ref{b1}) can be found by solving this linear
third-order ODE. Finally, these solutions can be used as a basis
for constructing the linearizing transformations by solving the
quadratures
\begin{equation}\label{a14}
\begin{array}{ll}
\frac{\phi_x}{\phi}=z-z_1,&\frac{\phi_y}{\phi}=w_1-w,\\
\frac{\psi_x}{\psi}=z-z_2,&\frac{\psi_y}{\psi}=w_2-w.\\
\end{array}
\end{equation}

The aim of this paper is to investigate the linearization of
second order ODEs using $\lambda$-symmetries. The outline of the
paper is as follows. The alternative proof of Lie's approach for
linearization of second order ODEs is provided in Section 2. This
section also presents a new $\lambda$-symmetry linearization
criteria which provides an alternative approach for constructing
the linearization transformations. The relationship between
$\lambda$-symmetries and the first integrals of ODEs play the key
role in the proving the results of Section 2. A familiarity with
standard results about the theory of $\lambda$-symmetries is
assumed and the reader is referred to the basic works
\cite{MR2001,MR2008,MR2009} on $\lambda$-symmetries. In Section 3,
we apply the new approach  to linearize the ODEs of the form
(\ref{a0}) with $F_3=F_2= 0$. Section 4 consists of examples
illustrating the application of the new approach to linearize ODEs
of the form (\ref{a0}) with $F_3\neq 0$ or $F_2\neq0$.
\section{Alternative proof of Lie's linearization approach $\&$ \\ new $\lambda$-symmetry linearization criteria}
 In this section, an alternative proof of Lie's approach to linearization of second order ODEs is presented. It is noticed that the relation between the
$\lambda$-symmetries
 and the first integrals provides a direct method to derive both
 the
 Lie's conditions and the quadratures. In addition, a $\lambda$-symmetry criteria for
linearization via a point transformation is stated. This criteria
provides a new approach for constructing the linearizing transformations, utilizing $\lambda$-symmetries.

The two first integrals
\begin{equation}\label{a5}
\begin{array}{ll}
 I_1=u_t,&I_2=u-t~u_t\\
\end{array}
\end{equation}
of the free particle equation (\ref{a1})  take the form
\begin{equation}\label{a6}
\begin{array}{ll}
 I_1={\frac {\psi_x  +\psi_y   y'}{\phi_x  +  \phi_y   y'}},&I_2=\psi-\phi \left({\frac {\psi_x  +\psi_y   y'}{\phi_x  +  \phi_y   y'}}\right),\\
\end{array}
\end{equation}
when expressed in the new variables $x$ and $y(x)$ defined by
equation  (\ref{a2}). It follows that all linearizable equations
(\ref{a0}) should have the two first integrals (\ref{a6}). Hence,
 the relationship between the first integrals and the
$\lambda$-symmetries \cite{MR2008,MR2009} implies that  all
linearizable equations (\ref{a0}) should have the two
$\lambda$-symmetries equivalent to the canonical pairs
$(\frac{\partial}{\partial y},\lambda_1)$ and
$(\frac{\partial}{\partial y},\lambda_2)$ given by
\begin{equation}\label{a7}
\begin{array}{ll}
\lambda_1=-\frac{{I_1}_y}{{I_1}_{y'}}=-A~{y'}^{2}- \left( B+w \right) y'-z\\
\lambda_2=-\frac{{I_2}_y}{{I_2}_{y'}}=-A~{y'}^{2}- \left( B+w-\frac{\phi_y}{\phi} \right) y'-z+\frac{\phi_x}{\phi}, \\
\end{array}
\end{equation}
where $A, B, w$ and $z$ are given by equation (\ref{a4})

The free particle equation (\ref{a1}) possesses the functionally dependent quotient first
integral $I_3=\frac{I_2}{I_1}=\frac{u-t~u_t}{u_t}$ as well. It is
worthwhile to mention here that the three triplets of Lie algebras
of $I_1$, $I_2$ and $I_3$ which are isomorphic to each other
generate the Lie algebra $sl(3,R)$ of the free particle equation
\cite{Lea1}. Similar to the above, it can be seen that all
linearizable equations (\ref{a0}) should have the third associated
$\lambda$-symmetry equivalent to the canonical pair
$(\frac{\partial}{\partial y},\lambda_3)$, where
\begin{equation}\label{a8}
\begin{array}{ll}
\lambda_3=-\frac{{I_3}_y}{{I_3}_{y'}}=-A~{y'}^{2}- \left( B+w-\frac{\psi_y}{\psi} \right) y'-z+\frac{\psi_x}{\psi}. \\
\end{array}
\end{equation}
The expressions for  $\lambda_1,\ \lambda_2$ and $\lambda_3$  can
be simplified, using the system (\ref{a3}), as
\begin{equation}\label{a9}
\begin{array}{ll}
\lambda_1=-F_3~{y'}^{2}- \left( F_2-w \right) y'-z\\
\lambda_2=-F_3~{y'}^{2}- \left( F_2-w-\frac{\phi_y}{\phi} \right) y'-z+\frac{\phi_x}{\phi} \\
\lambda_3=-F_3~{y'}^{2}- \left( F_2-w-\frac{\psi_y}{\psi} \right) y'-z+\frac{\psi_x}{\psi}. \\
\end{array}
\end{equation}

Employing the definition of $\lambda$-symmetry \cite{MR2001}, the
canonical pair $(\frac{\partial}{\partial y},\lambda_1)$ leads to
the following system for $w$ and $z$
\begin{equation}\label{a10}
\begin{array}{ll}
w_y=-  w^{2} +  {\it F_2} w +  {\it F_3}z+ {\frac {\partial F_3}{\partial x}}  -{\it F_1} {\it F_3},  \\
w_x-z_y=2\,w  z -2\,{\it F_3}  F  +{\frac {\partial F_2}{\partial x}}   -{\frac {\partial F_1}{\partial y _{{}}}}\\
z_x=z^{2} -  {\it F_1} z -  F w  + {\frac {\partial F}{\partial y_{{}}}} +F {\it F_2}.  \\
\end{array}
\end{equation}
Similarly, it follows by applying the definition of
$\lambda$-symmetry for $\lambda_2$ and $ \lambda_3$, along with
using the system (\ref{a10}), that the transformations $\phi$ and
$\psi$ can be given by finding two linearly independent non-constant
solutions of the following over-determined system
\begin{equation}\label{a11}
\begin{array}{ll}
S_{yy}  +\left( 2\,w-{\it F_2} \right) S_y +{\it F_3}\,S_x =0\\
S_{xy}  + w~S_x -z~S_y  =0\\
S_{xx}  +\left( {\it F_1}-2\,z \right) S_x -F~S_y  =0\\
\end{array}
\end{equation}
Based on the above discussion we have obtained the following two
$\lambda$-symmetry criteria for linearization via point
transformations.
\begin{theorem}
A scalar second-order ODE (\ref{a}) is linearizable via point
transformations (\ref{a2}) if and only if  it has the cubic in
first derivative form (\ref{a0})
with the $\lambda$-symmetries equivalent to the canonical pair $(\frac{\partial}{\partial y},\lambda_1)$   for $\lambda_1=-F_3~{y'}^{2}- \left( F_2-w \right) y'-z$
and the transformations $\phi$ and $\psi$ satisfying equation
(\ref{a11}), where $w$ and $z$ are auxiliary functions.
\end{theorem}
\proof The proof in one direction is given in the preceding discussion, so we  prove the other way.

Assume that equation (\ref{a0}) admits the $\lambda$-symmetries
equivalent to the canonical pair $(\frac{\partial}{\partial
y},\lambda_1)$   for $\lambda_1=-F_3~{y'}^{2}- \left( F_2-w
\right) y'-z$. Then, the system (\ref{a10}) for $w$ and $z$ is
obtained using the definition of $\lambda$-symmetry. Since the
transformations $\phi$ and $\psi$ satisfy equation (\ref{a11}),
the compatibility of the system (\ref{a11}), i.e. $S_{xy}=S_{yx}$,
$S_{xyx}=S_{xxy}$ and $S_{xyy}=S_{yyx}$, leads to the system
\begin{equation}\label{a12}
\begin{array}{l}
 \phi_y\left(-z_x+z^{2} -  {\it F_1} z -  F w  + {\frac {\partial F}{\partial y_{{}}}} +F {\it F_2},\right)-\phi_x  \left(-2\,z_y  -z w-w_x+ {\frac {\partial F_1}{\partial y}} +F  {\it F_3}     \right) =0\\
 \phi_x  \left(w_y+  w^{2} -  {\it F_2} w -  {\it F_3}z- {\frac {\partial F_3}{\partial x}}  +{\it F_1} {\it F_3}\right) -\phi_y  \left( 2\,w_x  -z  w+z_y   -{\frac {\partial F_2}{\partial x}}   +F  {\it F_3} \right)=0.\\
\end{array}
\end{equation}
Substituting the system (\ref{a10}) into the system (\ref{a12})
and noting that the Jacobian \\$ J=\phi_x\psi_y   - \phi_y
\psi_x\neq0 $, one finds
\begin{equation}\label{a13}
w_x+z_y=\frac{1}{3} \left({\frac {\partial F_2}{\partial x}}+{\frac {\partial F_1}{\partial y}}\right). \\
\end{equation}
Finally, system (\ref{a10}) and equation (\ref{a13}) are
equivalent to the Lie's conditions (\ref{b1}) for linearizable
equations and so the compatibility, $w_{xy}=w_{yx}$  and
$z_{xy}=z_{yx}$, of the Lie's conditions (\ref{b1}) leads to the
invariant equations (\ref{b2}).
\endproof
\begin{corollary}
A scalar second-order ODE (\ref{a}) is linearizable via point
transformations if and only if  it has the cubic in first
derivative form (\ref{a0}) with the $\lambda$-symmetries
equivalent to the canonical pair $(\frac{\partial}{\partial
y},\lambda_1)$ for $\lambda_1=-F_3~{y'}^{2}- \left( F_2-w \right)
y'-z$, where $w$ and $z$ are auxiliary functions satisfying the
equation
$w_x+z_y=\frac{1}{3} \left({\frac {\partial F_2}{\partial x}}+{\frac {\partial F_1}{\partial y}}\right)$.\\
\end{corollary}
\begin{remark}\rm
The compatibility of the system (\ref{a11}) is guaranteed by both
 the system (\ref{a10}) and equation
(\ref{a13}), i.e. $S_{xy}=S_{yx}$, $S_{xyx}=S_{xxy}$ and $S_{xyy}=S_{yyx}$.
Hence for each $w$ and $z$ given by solving (\ref{a10}) and
(\ref{a13}), one can construct transformations by solving the
over-determined system (\ref{a11}). This provides a new approach
for constructing the linearizing transformations whose
implementation  requires finding only one solution of the system
(\ref{b1}). In comparison, the standard implementation of Lie's
linearization approach involves solving Lie's quadratures
(\ref{a14}) which requires finding three solutions of the
system (\ref{b1}).
\end{remark}
Finally, in order to obtain an alternative proof of Lie's
quadratures (\ref{a14}) using $\lambda$-symmetries, it is noticed
that $\lambda_1,\ \lambda_2$ and $\lambda_3$ given by (\ref{a9})
can be written as
\begin{equation}
\begin{array}{ll}
\lambda_1=-F_3~{y'}^{2}- \left( F_2-w \right) y'-z,\\
\lambda_2=-F_3~{y'}^{2}- \left( F_2-w_1 \right) y'-z_1,\\
\lambda_3=-F_3~{y'}^{2}- \left( F_2-w_2 \right) y'-z_2,\\
\end{array}
\end{equation}
where $ w_1=w+\frac{\phi_y}{\phi},
z_1=z-\frac{\phi_x}{\phi},w_2=w+\frac{\psi_y}{\psi}$ and
$z_2=z-\frac{\psi_x}{\psi}$. Therefore, using the definition of
$\lambda$-symmetry shows that $(z, w), (z_1, w_1)$ and $(z_2,w_2)$
are the three solutions for the system (\ref{b1}). This completes
the alternative proof of Lie's approach.
\section{Linearization of the ODE's of the form (\ref{a0}) with
$F_3=F_2= 0$} The non-linear second order ODEs of the form
\begin{equation}\label{form1}
y''+F_1(x,y)y'+F(x,y)=0
\end{equation}
satisfy the Lie-Tress\'e linearization criteria (\ref{b2}) if and only
if
$$
F_1(x,y) =  a ( x  ) y  +b ( x ), \quad F(x,y)= \frac{ \left( a (
x )  \right) ^{2}}{9}
  y^{3}+ \frac{1}{3} \left( {\frac {da(x)}{dx}} +a ( x ) b ( x )  \right)
  y ^{2}+c( x ) y + d(x)
$$
Here we consider the class of nonlinear equations
\begin{equation}\label{main}
{\frac {d^{2}y}{d{x}^{2}}} + \left( a \left( x
 \right) y  +b \left( x \right)  \right) {\frac {dy}{dx
}} +\frac{1}{9} \left( a \left( x \right)  \right) ^{2}
  y^{3}+ \frac{1}{3} \left( {\frac {da(x)}{dx}} +a \left( x \right) b \left( x \right)  \right)
  y ^{2}+c \left( x \right) y +d(x)=0
\end{equation}
where $a(x)\neq0$.

Solving the system (\ref{b1}) gives  the three solutions $(w,z),
(w_1,z_1)$ and $(w_2,z_2)$
\begin{equation}\label{raed}
\begin{array}{ll}
w(x,y)={\frac {a \left( x \right) g_1\left( x \right) }{ ya \left( x \right) g_1 \left( x \right)-3\,{g_1}' \left( x\right) }},&z(x,y)=\frac{1}{3}\,{\frac {9\, a \left( x \right){g_1}''\left( x \right)   - \left( 9\,a' \left( x \right) +6\, y{ a \left( x \right)}^{2}\right) {g_1}' \left( x \right) + {y}^{2}{ a \left( x\right)  } ^{3}g_1 \left( x \right) }{a \left( x \right)\left( ya \left( x \right) g_1 \left( x \right) -3\,{g_1}' \left( x \right)  \right) }},\\
w_1(x,y)={\frac {a \left( x \right) g_2\left( x \right) }{ ya \left( x \right) g_2 \left( x \right)-3\,{g_2}' \left( x\right) }},&z_1(x,y)=\frac{1}{3}\,{\frac {9\, a \left( x \right){g_2}''\left( x \right)   - \left( 9\,a' \left( x \right) +6\, y{ a \left( x \right)}^{2}\right) {g_2}' \left( x \right) + {y}^{2}{ a \left( x\right)  } ^{3}g_2 \left( x \right) }{a \left( x \right)\left( ya \left( x \right) g_2 \left( x \right) -3\,{g_2}' \left( x \right)  \right) }},\\
w_2(x,y)={\frac {a \left( x \right) g_3\left( x \right) }{ ya \left( x \right) g_3 \left( x \right)-3\,{g_3}' \left( x\right) }},&z_2(x,y)=\frac{1}{3}\,{\frac {9\, a \left( x \right){g_3}''\left( x \right)   - \left( 9\,a' \left( x \right) +6\, y{ a \left( x \right)}^{2}\right) {g_3}' \left( x \right) + {y}^{2}{ a \left( x\right)  } ^{3}g_3 \left( x \right) }{a \left( x \right)\left( ya \left( x \right) g_3 \left( x \right) -3\,{g_3}' \left( x \right)  \right) }},\\
\end{array}
\end{equation}
where $g_i(x), i=1,2,3,$ are the three linearly independent
solutions of the linear third order ODE
\begin{equation}\label{MODE}
\begin{array}{ll}
Y'''(x)+\left(b(x)-2\frac{a'(x)}{a(x)}\right)Y''(x)-\left(\frac{a''(x)}{a(x)}-2{\frac{a'(x)}{a(x)}}^2+b(x)\frac{a'(x)}{a(x)}-c(x)\right)Y'(x)+\frac{d(x)a(x)}{3}Y(x)=0.\\
\end{array}
\end{equation}
Now, in order to construct the linearizing transformation using
Lie's linearization approach, one should solve the quadratures
(\ref{a14}) by considering the three solutions $(w,z),$
$(w_1,z_1)$ and $(w_2,z_2)$. However, our approach requires
utilizing any one the three solutions $(w,z),$ $(w_1,z_1)$,
$(w_2,z_2)$ in order to determine two non-constant solutions of
the over-determined system (\ref{a11}) which yield the
linearization transformations.

As an illustration, we consider ODEs (\ref{main}) with $d(x)=0$.
For this case $g_1(x)=1$ is a constant solution of the ODE
(\ref{MODE}). So equation (\ref{raed}) gives the solution
$(w,z)=\left(\frac{1}{y},\frac{1}{3}ya(x)\right)$ of the system
(\ref{b1}). Now, solving the over-determined system (\ref{a11})
gives the linearizing point transformations
\begin{equation}\label{MTR}
\begin{array}{ll}
\phi(x,y)=\frac{1}{3}\int a(x) h_1(x) dx -\frac{h_1(x)}{y},&\psi(x,y)=\frac{1}{3}\int a(x) h_2(x) dx -\frac{h_2(x)}{y}.\\
\end{array}
\end{equation}
where $h_i(x), i=1,2,$ are the two linearly independent solutions
of the linear second order ODE
\begin{equation}\label{MODE2}
\begin{array}{ll}
H''(x)+b(x)H'(x)+c(x)H(x)=0.\\
\end{array}
\end{equation}
Finally, since $u(t)=c_1t+c_2$ is the general solution of the free
particle equation $u_{tt}=0$,  the general solution of the ODE
(\ref{main}) can be given as
\begin{equation}\label{MGS}
\frac{1}{3}\int a(x) h_2(x) dx -\frac{h_2(x)}{y}=c_1\left(\frac{1}{3}\int a(x) h_1(x) dx -\frac{h_1(x)}{y}\right)+c_2.\\
\end{equation}
\begin{examples}\rm
We consider the ODE (\ref{main}) with $a(x)=3,b(x)=0$ and $c(x)=0$
that gives the modified Emden equation
\begin{equation}
y''+3yy'+y^3=0,
\end{equation}
for which the ODE (\ref{MODE2}) reduces to
\begin{equation}
H''(x)=0.\\
\end{equation}
Hence the linearizing point transformations can be written using
(\ref{MTR}) as
\begin{equation}
\begin{array}{ll}
\phi(x,y)=x-\frac{1}{y},&\psi(x,y)=\frac{x^2}{2}-\frac{x}{y}.\\
\end{array}
\end{equation}
Finally the general solution can be written using (\ref{MGS}) as
\begin{equation}
y(x)=\frac{2x+c_1}{x^2+c_1 x+c_2}.\\
\end{equation}
\end{examples}
\begin{examples}\rm
For  $a(x)=3k,b(x)=b$ and $c(x)=\frac{b^2}{4}$, ODE (\ref{main})
gives the $\rm{Li\acute{e}nard}$ type equation
\begin{equation}
y''+(b+3ky)y'+k^2y^3+bky^2+\frac{b^2}{4}y=0,
\end{equation}
and the ODE (\ref{MODE2}) reduces to
\begin{equation}
H''(x)+bH'(x)+\frac{b^2}{4}H(x)=0.\\
\end{equation}
So the linearizing point transformations can be stated, using
(\ref{MTR}), as
\begin{equation}
\begin{array}{ll}
\phi(x,y)=\frac{2ky+b}{bky}e^{-\frac{b}{2}x},&\psi(x,y)=\frac{2bkxy+4ky+b^2x}{b^2ky}e^{-\frac{b}{2}x},\\
\end{array}
\end{equation}
which lead to the general solution via (\ref{MGS}) as
\begin{equation}
y(x)=\frac{b^2(c_1-x)}{2bkx+4k-2c_1bk+c_2b^2k e^{\frac{b}{2}x}}.
\end{equation}
\end{examples}
\section{Examples of linearization of the ODEs of the form (\ref{a0})
with $F_3 \ne 0$ or $F_2 \ne 0$ } In this section, we illustrate
the application of our approach to obtain linearization
transformations for ODEs that are cubic or quadratic in the first
derivative.
\begin{examples}\rm
As the first example, we consider the nonlinear ODE
\begin{equation}\label{ex41}
y''-\frac{2}{x+y}{y'}^2-\frac{1}{x+y}y'=0
\end{equation}
which satisfies the  Lie-Tress\'e linearization test (\ref{b2}). The
three solutions $(w,z), (w_1,z_1)$ and $(w_2,z_2)$ of the Lie's
conditions (\ref{b1}) are
\begin{equation}
\begin{array}{ll}
w(x,y)=-\frac{1}{x+y},& z(x,y)=0,\\
w_1(x,y)=\frac{x}{y(x+y)},& z_1(x,y)=0,\\
w_2(x,y)=\frac{x}{(x+2y)(x+y)},& z_2(x,y)=-\frac{2(x+y)}{x(x+2y)}.\\
\end{array}
\end{equation}
Our approach requires only one solution of the Lie's conditions
(\ref{b1}). Solving the over-determined system (\ref{a11}), by
considering only $(w,z)$, gives the two non-constant solutions
which results in the point transformations
\begin{equation}
\begin{array}{ll}
\phi(x,y)=y,&\psi(x,y)=x(x+2y).\\
\end{array}
\end{equation}
that linearize the ODE (\ref{ex41}) to the free particle equation
$u_{tt}=0$. So the general solution of the ODE (\ref{ex41}) is
\begin{equation}
x(x+2y)=c_1y+c_2.\\
\end{equation}
It is worth mentioning that any of the other solutions $(w_1,z_1)$
and $(w_2,z_2)$ of the Lie's conditions (\ref{b1}) can be used to
derive different linearizing point transformations that will lead
to the same general solution.
\end{examples}
\begin{examples}\rm
The nonlinear ODE
\begin{equation}\label{ex42}
xy''-{y'}^3-y'=0
\end{equation}
satisfies the  Lie-Tress\'e linearization test (\ref{b2}). A solution
$(w,z)$ of the Lie's conditions (\ref{b1}) is
\begin{equation}
\begin{array}{ll}
w(x,y)=\frac{1}{y},& z(x,y)=0.\\
\end{array}
\end{equation}
Now, solving the over-determined system (\ref{a11}) gives the two
non-constant solutions which gives rise to the point transformation
\begin{equation}
\begin{array}{ll}
\phi(x,y)=\frac{1}{y},&\psi(x,y)=y+\frac{x^2}{y}\\
\end{array}
\end{equation}
that linearize the ODE (\ref{ex42}) to the free particle equation
$u_{tt}=0$. So the general solution of the ODE (\ref{ex42}) is
obtained as
\begin{equation}
x^2 + y^2=c_1 +c_2 y.\\
\end{equation}
\end{examples}
\begin{examples}\rm
In this example, we find the general solution of the nonlinear ODE
\begin{equation}\label{ex42}
y''-\frac{1}{x}{y'}^3+\frac{2y}{y^2-1}{y'}^2-\frac{1}{x}y'=0
\end{equation}
via linearization transformations. ODE (\ref{ex42}) satisfies the
Lie-Tress\'e linearization test (\ref{b2}) and a solution $(w,z)$ of the
Lie's conditions (\ref{b1}) is
\begin{equation}
\begin{array}{ll}
w(x,y)=\frac{3y^2-3}{y^3-3y},& z(x,y)=0.\\
\end{array}
\end{equation}
Now, solving the over-determined system (\ref{a11}) gives the two
non-constant solutions which gives the point transformation
\begin{equation}
\begin{array}{ll}
\phi(x,y)=\frac{1}{y(y^2-3)},&\psi(x,y)=\frac{(y^3-3y+2)\ln(y-1) + (2+3y-y^3)\ln(y+1) +3 x^2-y^2   }{y(y^2-3)}\\
\end{array}
\end{equation}
that linearize the ODE (\ref{ex42}) to the free particle equation
$u_{tt}=0$. Hence the general solution of the ODE (\ref{ex42}) can be
written as
\begin{equation}
(y^3-3y+2)\ln(y-1) + (2+3y-y^3)\ln(y+1) +3 x^2-y^2   =c_1  +c_2 y(y^2-3) .\\
\end{equation}
\end{examples}
\section{Conclusion}
The question of linearization of second order ODEs is investigated
employing $\lambda$-symmetries. The relationship between
$\lambda$-symmetries and the first integrals plays an important
role and provides a direct method to derive both of
 Lie's conditions (\ref{b1}) and quadratures (\ref{a14}). The relationship
 further leads to a $\lambda$-symmetry criteria for linearization via
point transformations. This criteria provides a new approach for
constructing the linearizing transformations whose implementation
requires finding only one solution of the Lie's conditions
(\ref{b1}). In comparison, the standard implementation of Lie's
linearization approach involves solving Lie's quadratures
(\ref{a14}) which requires finding three solutions of the Lie's
conditions (\ref{b1}).

It is expected that the relationship between $\lambda$-symmetries
and the first integrals of ODEs can play a significant role in
deriving new linearization criteria for higher order ODEs.
\subsection*{Acknowledgments}
 The authors would like to thank the King Fahd University of
Petroleum and Minerals for its support and excellent research
facilities.

\end{document}